\documentclass[a4paper,12pt,leqno]{article}
\usepackage{latexsym}
\usepackage[all]{xy}

\usepackage{amssymb} 
\usepackage{amsmath} 

\usepackage{tikz}
\usetikzlibrary{arrows,cd}  
\usepackage{amsthm}

\def\P{{\mathbf{P}}}
\def\Z{{\mathbb{Z}}}
\def\K{{\mathbb{K}}}
\def\CC{{\mathbb{C}}}
\def\R{{\mathbb{R}}}
\def\A{{\mathcal{A}}}
\def\B{{\mathcal{B}}}

\DeclareMathOperator{\codim}{codim}

\DeclareMathOperator{\Der}{Der}

\DeclareMathOperator{\pd}{pd}
\DeclareMathOperator{\depth}{depth}

\DeclareMathOperator{\res}{res}

\numberwithin{equation}{section}

\newcommand{\owari}{\hfill$\square$}

\newtheorem{theorem}{Theorem}[section]
\newtheorem{prop}[theorem]{Proposition}
\newtheorem{cor}[theorem]{Corollary}
\newtheorem{lemma}[theorem]{Lemma}
\newtheorem{define}[theorem]{Definition}

\theoremstyle{remark}
\newtheorem{rem}[theorem]{Remark}
\newtheorem{example}[theorem]{Example}

\title{Tame arrangements}
\author{Takuro Abe}
\date{\today}

\begin{document}

\maketitle

\begin{abstract}
Tame arrangements were informally introduced by Orlik and Terao for the study of Milnor fibers of hyperplane arrangements. After that, tame arrangements have been applied to a lot of researches on arrangements including freeness, master functions and critical varieties, 
Solomon-Terao algebras, D-modules, Bernstein-Sato polynomials and 
likelihood geometry. Though arrangements are generically tame, the research on tame arrangements themselves have been only few. In this article we establish foundations for the  research of tame arrangements. Namely, we prove the addition theorem for tame arrangements, Ziegler-Yoshinaga type results for tameness and combinatorially determined  tameness. \end{abstract}

\section{Introduction}
Let $\K$ be a field, $V=\K^\ell$ and $\A$ be a hyperplane arrangement in $V$, i.e., 
a finite set of linear hyperplanes in $V$. For $H \in \A$ let us fix a defining linear form $\alpha_H \in V^*$. Let $S=\K[x_1,\ldots,x_\ell]$ be a coordinate ring of $V$, 
$\Der S:=\oplus_{i=1}^\ell S\partial_{x_i}$ be the derivation module, 
and let $
\Omega^1_V:=\oplus_{i=1}^\ell S dx_i
$ be a free module of regular differential $1$-forms. Then we can define the regular $p$-forms as follows:
$$
\Omega^p_V:=\wedge^p \Omega^1_V.
$$
Let $Q(\A):=\prod_{H \in \A} \alpha_H$. Now we can define the main object in this article:

\begin{define}
Let $0 \le p \le \ell$. The module of logarithmic $p$-forms $\Omega^p(\A)$ of $\A$ is defined by 
$$
\Omega^p(\A):=\{\omega 
\in \frac{1}{Q(\A)}\Omega^p_V \mid 
Q(\A)d\alpha_H \wedge \omega \in \alpha_H \Omega^{p+1}_V\ 
(\forall H \in \A)\}.
$$
Note that $\Omega^0(\A)=S$ and $\Omega^\ell(\A) =
S\frac{1}{Q(\A)} \wedge_{i=1}^\ell dx_i$.
\label{logformdefine}
\end{define}

In general $\Omega^p(\A)$ is not a free but reflexive $S$-graded module of rank $\frac{\ell!}{p!(\ell-p)!}$, see \cite{ST} or \cite{OT} for example. So we say that $\A$ is \textbf{free} if $\Omega^1(\A)$ is a free module of rank $\ell$. If $\A$ is free, then it is known that $\Omega^p(\A)$ coincides with $\wedge^p \Omega^1(\A)$ (Theorem \ref{MS}), so $\Omega^p(\A)$ is free too (Proposition \ref{freewedge}). Free arrangements have been in the center of the algebraic research of hyperplane arrangements in this 40 years.
On the other hand, recently, tame arrangements are also playing central roles in several studies on arrangements like algebra, geometry and analysis:

\begin{define}
$\A$ is \textbf{tame} if the projective dimension $\pd_S \Omega^p(\A)$ of $\Omega^p(\A)$ as an $S$-module is at most $p$ for all $0\le p\le \ell$.
\label{tamedefinition}
\end{define}

Tame arrangements were originally defined by Orlik and Terao in \cite{OT2} to investigate the critial points of the Milnor fiber of hyperplane arrangements. After that, there have been several researches using tameness. For example, Schulze \cite{Sch} and the author \cite{A0} for the freeness criterion of tame arrangements, Denham-Schulze for local freeness \cite{DSch}, Cohen-Denham-Falk-Varchenko \cite{CDFV}, 
Denham-Schenck-Schulze-Wakefield-Walther \cite{DSSWW}, 
Denahm-Garrousian-Schulze \cite{DGS} and Denham-Steiner \cite{DSt} for the master functions, logarithmic ideals and critical varieties, Maeno-Murai-Numata and the author \cite{AMMN} for the Solomon-Terao polynomials and algebras, 
and Walther \cite{W} and Bath \cite{B1}, \cite{B2} for D-modules and the Bernstein-Sato polynomials of hyperplane arrangements. The range where tame arrangements appear is increasing, for example, a recent preprint \cite{KKMSW} by 
Kahle, K\"{u}hne, M\"{u}hlherr, Sturmfels and Wiesmann, 
shows that tameness is also related to likelihood geometry.
We can see from these results that tameness is a very useful and widely-used property of arrangements. In fact, it is known that tameness is a generic property, i.e., both generic arrangements and free arrangements are known to be tame (\cite{RT}, \cite{ST}).

However, there have been few researches on tame arrangements themselves. It is easy to show that every arrangement in $\K^\ell$ for $\ell \le 3$ is tame by reflexivity (Theorem \ref{dualfree}). Except for it, the only result to determine the tameness is the following wonderful result due to Musta\c{t}\u{a} and Schenck. Let us recall it after defining local freeness of arrangements.

\begin{define}
(1)\,\,Let $X \subset V$. We say that $\A$ is 
\textbf{locally free along $X$} if $\A_x:=\{H \in \A \mid x \in H\}$ is free for all $0 \neq x \in X$.

(2)\,\, We say that $\A$ is \textbf{locally free} if $\A$ is locally free along $V \setminus \{0\}$.
\label{locallyfree}
\end{define}

\begin{theorem}[Theorem 5.3, \cite{MS}]
Assume that $\A$ is locally free and $\pd_S \Omega^1(\A)=1$. Then 
$\A$ is tame. Explicitly, it holds that $
\wedge^p \Omega^1(\A)=\Omega^p(\A),\ \pd_S \Omega^p(\A)=p$ for $0 \le p \le \ell-2$ and 
$\pd_S \Omega^{\ell-1}(\A)=\ell-2$. 
\label{MS}
\end{theorem}

To show Theorem \ref{MS}, they used the sheaf isomorphism and the explicit free resolution of the wedge product in this case called the Lebelt resolution, by which they could give explicit minimal free resolutions. 

The aim of this article is to construct a foundation for the theory of tame arrangements. First we introduce an inductive way to construct a tame arrangement from a given tame arrangement, so this is an addition-deletion type approach whose origin is Terao's addition-deletion theorem for free arrangements \cite{T1}. The first main result in this article is the following addition theorem for tameness:

\begin{theorem}[Addition theorem for tameness I]
Let $\A=\A' \cup \{H\}$. If $\A'$ is tame and both $\A'$ and 
$\A$ are locally free 
along $H$, then $\A$ and $\A^H$ are both tame.
\label{addtame}
\end{theorem}

Theorem \ref{addtame} is interesting since a local freeness along $H$ and the tameness of $\A \setminus \{H\}$ confirm the tameness of both $\A$ and $\A^H$. 
Also, we can show more addition-theorem type assertion as follows:

\begin{theorem}[Addition theorem for tameness II]
Assume that $\A^H$ is tame, and 
$\A'$ is locally free along $H$. Then $\A $ is tame if 
$\A':=\A \setminus \{H\}$ is tame.
\label{adddeltame}
\end{theorem}

Tameness has been a very useful but hard property arrangement to check. However, by Theorems \ref{addtame} and \ref{adddeltame}, now we can look for and construct a lot of tame arrangements easily. Let us show some  examples of tame arrangements made by these theorems.

\begin{example}
(1)\,\,
We say that $\A$ is \textbf{generic} if 
for any $H_1,\ldots,H_k \in \A$ with $ k \le \ell-1$, it holds that 
$\codim_V \bigcap_{i=1}^kH_i =k$. 
It was known that a generic arrangement is tame, which is based on the explicit construction of generators and minimal free resolutions due to \cite{RT}. 
Now we can show it by using Theorem \ref{addtame} easily. By definition, we can construct any generic arrangement by adding a generic hyperplane to 
a generic arrangement. So the first step to make it is to add a generic hyperplane to a free arrangement, which is clearly free by, e.g., Proposition \ref{freewedge}. So let us assume that $\A'$ is generic, thus tame and let $H$ be generic with respect to $\A'$. 
Let $\A:=\A' \cup\{H\}$. Since $\A'$ and $\A$ are generic, clearly they are both locally free along $H$. So Theorem \ref{addtame} shows that both 
$\A$ and $\A^H$ are tame.

(2)\,\,
Let $\A'$ be defined by 
$$
Q(\A')=x_{\ell+1}(x_1+\cdots+x_\ell)\prod_{i=1}^\ell x_i
$$
in $\R^{\ell+1}$. $\A'$ is tame since it is a coning of a generic, thus a tame arrangement $(x_1+\cdots+x_\ell)\prod_{i=1}^\ell x_i=0$. Note that $\A'$ is not locally free since $\A'$ contains that generic arrangement as its localization. 
Let us add $H:x_1+x_{\ell+1}=0$ to $\A'$ to obtain $\A:=\A' \cup \{H\}$. It is easy to show that $\A'$ is locally free along $H$. Also, $\A^H$ is tame since $\A^H$ is also generic. So Theorem \ref{adddeltame} shows that $\A$ is tame which is not locally free. 

(3)\,\,
Let $\A$ be defined by 
$$
xyzwu(x+y+z+w+u)(x+y+z)=0
$$
in $\CC^5$. Let $\A \ni H:x+y+z=0$ and let $\A':=\A \setminus \{H\}$.
Since $\A'$ is generic, it is tame. Also, since $\A^H$ is defined by 
$$
yz(y+z)wu(w+u)=0
$$
which is a direct sum of two-dimensional (thus free) arrangements, it is free, in particular, it is tame. Since a generic arrangement $\A'$ is locally free, Theorem \ref{adddeltame} shows that $\A$ is tame. Note that 
$\A$ is not free since for $X=\{x=y=z=0\} \in L(\A)$, $\A_X$ is generic thus not free.
\end{example}

Moreover, we can show the restriction theorem for tameness as follows.

\begin{theorem}[Restriction theorem for tameness]
If $\A$ and $\A'$ are both tame, and $\A$ is locally free along $H$, then 
$\A^H$ is tame.
\label{resttame}
\end{theorem}

Also as in \cite{Z}, \cite{Y1} and \cite{Y2}, nowadays it is usual to study the relation of several properties of logarithmic differential forms of the original arrangement and its Ziegler restriction. We show the following Ziegler-type and Yoshinaga-type results on tameness, which are the second main results in this article (see Definitions \ref{multiarrangement} and \ref{Zrest} for the definitions appearing below):

\begin{theorem}
Let $H \in \A$. 
Assume that $\A$ is tame and $\A$ is locally free along $H$. Then the Ziegler restriction $(\A^H,m^H)$ is tame, i.e., 
$\pd_{\overline{S}} \Omega^p(\A^H,m^H) \le p$ for $0 \le p \le \ell-1$, where $\overline{S}:=S/\alpha_H S$.
\label{Zresttame}
\end{theorem}

\begin{theorem}
\begin{itemize}
    \item[(1)]
    Let $\ell=4$. Then $\A$ is tame if there is some $H \in \A$ such that 
    $\A \setminus \{H\}$ is tame, and locally free along $H$.
\item[(2)]
Let $\ell \ge 5$, and $H \in \A$. Assume that the Ziegler restriction $(\A^H,m^H)$ is tame and $\A$ is locally free along $H$.
Then $\A$ is tame.
\end{itemize}
\label{Ytame}
\end{theorem}

Note that the statements are trivial when $\ell \le 3$ since all such multiarrangements are tame. Also, because of that, the tameness of the Ziegler restriction of arrangements in $\K^4$ is automatical. So we need to divide the statement depending on $\ell =4$ or $\ell \ge 5$, which is similar to the original Yoshinaga's criterion which is also divided into two statamenets depending on whether $\ell =3$ or $\ell \ge 4$. 

When $\ell \ge 4$, we can regard Yoshinaga's criterion (Theorem \ref{Ycriterion}) as follows. Assume that $\pd_{\overline{S} }\Omega^1(\A^H,m^H)=0$ and $\A$ is locally free along $H$, then $\pd_S \Omega^1(\A)=0$. So Theorems \ref{Zresttame} and \ref{Ytame} can be regarded as a generalized Yoshinaga's criterion as follows.

\begin{cor}
Assume that $\ell \ge 5$ and $\A$ is locally free along $H$. Then 
$\A$ is tame 
if and only if $(\A^H,m^H)$ is tame. Explicitly, the following are equivalent:
\begin{itemize}
\item[(1)]
$\A$ is tame.
\item[(2)]
$\pd_S \Omega^1(\A) \le 1$.
\item[(3)]
$(\A^H,m^H)$ is tame.
\item[(4)]
$\pd_{\overline{S} } \Omega^1(\A^H,m^H) \le 1$.
\end{itemize}
\label{generalizedY}
\end{cor}

We give one example of tame arrangements by using Theorem \ref{Ytame} and 
Theorem \ref{multigenerictame} in \S4.

\begin{example}
Let $\A$ be an irreducible and generic arrangement with $\ell \ge 5$, i.e., 
it holds that $\codim_V \bigcap_{i=1}^{\ell-1} H_i=\ell-1$ for any distinct hyperplanes $H_1,\ldots,H_{\ell-1} \in \A$. Let us fix $\{z=0\} =H \in \A$ and for each $L \in \A$, fix a positive integer $n_L$, distinct elements $a_1^L,\ldots,a_{n_L}^L \in \K$, and define a new arrangement $\B$ by 
$$
\B:=H \cup \bigcup_{H \neq L \in \A } \bigcup_{i=1}^{n_L}\{\alpha_L=a_i^L z\}.
$$
Then $\B$ is no more generic, but we can show that $\B$ is generic. First, consider the Ziegler restriction $(\B^H,m^H)$ of $\B$ onto $H$, which is a generic multiarrangement. Theorem \ref{adddeltamemulti} shows that $(\B^H,m^H)$ is tame. Also, it is easy to show that $\A$ is locally free along $H$ because any localization at codimension $\ell-1$ is isomorphic to 
the free arrangement 
$$
z\prod_{i=1}^{\ell-1} \prod_j (x_i-a_j^Lz)=0.
$$
Hence Theorem \ref{Ytame} shows that $\A$ is tame.

Let us give more explicit example. Starting from a generic arrangement $\B$
$$
x_\ell(x_1+\cdots+x_\ell)\prod_{i=1}^{\ell-1} x_i=0$$
and deforming it into 
$$
\A:
x_{\ell}(x_1+\cdots+x_\ell)\prod_{i=1}^{\ell-2} \prod_{j=-10}^{10}x_i(x_i-(j-i)x_\ell)=0,
$$
we have no idea whether $\A$ is locally free, and so on. However, by the arguments above, we know that $\A$ is tame.
\end{example}

The organization of this article is as follows. In \S2 we recall several definitions and results for the proof of main results, and we give their proofs in \S3. In \S4 we give Ziegler-Yoshinaga type results for tameness, i.e., Theorems \ref{Zresttame} and \ref{Ytame}. 
\S5 is devoted for the investigation of the combinatorially determined tame arrangements so called inductively tame arrangements.
\medskip

\noindent
\textbf{Acknowledgements}.
The author 
is partially supported by JSPS KAKENHI Grant Numbers JP23K17298 and JP25H00399.

\section{Preliminaries}

From now on, for an arrangement $\A$, we assume that $$
\bigcap_{H \in \A} H=\{0\}.$$
Namely, all arrangements are essential unless otherwise specified. 
In this section let us introduce definitions 
and results necessary for the proofs of main results. 

\begin{define}
Let 
$$
L(\A):=\{\bigcap_{H \in \B} H \mid \B \subset \A\}
$$
be the \textbf{intersection lattice} and let 
$$
L_k(\A):=\{X \in L(\A)  \mid 
\codim_V =k\}
$$
for $ 0\le k \le \ell$.
\label{lattice}
\end{define}
Since we prove main theorems in a wider class, i.e., of multiarrangements, let us introduce several related terminologies. 

\begin{define}[\cite{Z}]
A function $m:\A\rightarrow \Z_{>0}$ is called the \textbf{multiplicity}, and the pair $(\A,m)$ is called the \textbf{multiarrangement}. Let $Q(\A,m):=\prod_{H \in \A} \alpha_H^{m(H)}$ and let $|m|:=\sum_{H \in \A} m(H)$. 

Let $1$ denote the constant multiplicity, i.e., $1(H)=1$ for all $H \in \A$. 
For two multiplicities $k,m$ on $\A$, $k \le m$ denotes that $k(H) \le m(H)$ for 
all $H \in \A$. Moreover, $k < m$ denotes that $k \le m$ and 
there is $H \in \A$ such that $k(H)<m(H)$. For $H \in \A$, let $\delta_H: \A 
\rightarrow \{0,1\}$ be the 
\textbf{characteristic multiplicity} defined by $\delta_H(L):=\delta_{H,L}$. 

We can define the 
\textbf{logarithmic differential $p$-forms} of $(\A,m)$ by 
\begin{eqnarray*}
\Omega^p(\A,m):= \{\omega \in \frac{1}{Q(\A,m)} \Omega^p_V 
\mid Q(\A,m)d\alpha_H \wedge \omega \in \alpha_H^{m(H)}\Omega^{p+1}_V\ 
(\forall H \in \A)\},
\end{eqnarray*}
and 
\textbf{logarithmic derivation module of order $p$} of $(\A,m)$ by 
\begin{eqnarray*}
D^p(\A,m):= \{\theta \in \wedge^p \Der S
\mid \theta(\alpha_H,f_2,\ldots,f_p)\in \alpha_H^{m(H)}S \ 
(\forall H \in \A,\ \forall f_i \in S)\}.
\end{eqnarray*}
\label{multiarrangement}
\end{define}

\begin{rem}
When $m\equiv 1$, $D^p(\A,1)=D^p(\A)$ and $\Omega^p(\A,1)=\Omega^p(\A)$. So if we define, or prove algebraic statements for multiarrangements, then it includes those for a usual arrangement $\A$.
\end{rem}

\begin{theorem}[\cite{Z}]
$D^p(\A,m)^* \simeq \Omega^p(\A,m)$ and $\Omega^p(\A,m)^* \simeq D^p(\A,m)$. Hence 
$\Omega^p(\A,m)$ and $D^p(\A,m)$ are both reflexive modules. In particular, $\pd_S D^p(\A) \le \ell-2$ and $\pd_S \Omega^p(\A) \le \ell-2$ for all $p$.
\label{dualfree}
\end{theorem}

\begin{define}
We say that $(\A,m)$ is \textbf{free} with exponents $\exp(\A,m)=(d_1,\ldots,d_\ell)$ if 
$D(\A,m) \simeq \bigoplus_{i=1}^\ell S[-d_i]$, equivalently, $\Omega^1(\A,m) \simeq 
\bigoplus_{i=1}^\ell S[d_i]$. 
\label{multifree}
\end{define}

\begin{prop}
If $\ell=2$, then $(\A,m)$ is free.
\label{2free}
\end{prop}

\noindent
\textbf{Proof}. Clear by the reflexitivity (Theorem \ref{dualfree}) and Auslander-Buchsbaum formula.\owari
\medskip

\begin{prop}[Lemma 1.3, \cite{ATW}]
If $\A$ is free, then $$
\wedge^p \Omega^1(\A,m) \simeq \Omega^p(\A,m)
$$
and 
$$
\wedge^p D(\A,m) \simeq D^p(\A,m).
$$
\label{freewedge}
\end{prop}

\begin{prop}[cf. Remark 2.4, \cite{AY}]
It holds that 
$$
D^p(\A,m) \simeq Q(\A,m) \Omega^{\ell-p}(\A,m)
$$
for $ 0 \le p \le \ell$.
\label{idenfify}
\end{prop}

\noindent
\textbf{Proof}. Apply the same proof as in Remark 2.4, \cite{AY} for multiarrangements.
\owari
\medskip

\begin{rem}
By Proposition \ref{idenfify}, several statements on $\Omega^p(\A,m)$ hold true for 
$D^p(\A,m)$.
\label{dualbothtrue}
\end{rem}

\begin{lemma}[Theorem 2.2, \cite{KS}]
Let $(\A,m)$ be a a multiarrangement and $x \subset S$ be a prime ideal. Then for $\A_x:=\{H \in \A \mid x \in H\}$ and $m_x:=m|_{\A_x}$, it holds that 
$$
\pd_S \Omega^p(\A_x,m_x) \le \pd_S \Omega^p(\A,m).
$$
In particular, for $ X \in L(\A)$ and $\A_X:=\{H \in \A \mid H \supset X\}$ and 
$m_X:=m|_{\A_X}$, it holds that 
$$
\pd_S \Omega^p(\A_X,m_X) \le \pd_S \Omega^p(\A,m).
$$
\label{pdlocalization}
\end{lemma}



For a multiarrangement, we can define its tameness too.

\begin{define}
We say that $(\A,m)$ is \textbf{tame} if 
$\pd_S \Omega^p(\A,m) \le p$ for all $p$.
\label{multitame}
\end{define}

Note that $\pd_S \Omega^p(\A,m) \le p$ holds for all $(\A,m)$ if $p=0,\ell-2,\ell-1,\ell$ by the definition of $\Omega^0(\A,m)$ and reflexivity Theorem \ref{dualfree}. 
Now let us introduce a restriction of a multiarrangement $(\A,m)$ onto $H \in \A$.

\begin{define}[Definition 0.2, \cite{ATW2}]
Let $(\A,m)$ be a multiarrangement, $H \in \A$.
Then we can define the \textbf{Euler restriction} 
$(\A^H,m^*)$ of $(\A,m)$ onto $H$ as follows. For $X \in \A^H$, let 
$m_X:=m|_{\A_X}$. Then $(\A_X,m_X)$ is free since $\A_X$ is essentially a two-dimensional arrangement, so let $\exp(\A_X,m_X)=(d_1^X,d_2^X,0\ldots,0)$. Then by \cite{ATW2}, we may assume that $\exp(\A_X,m_X-\delta_H)=(d_1^X-1,d_2^X,0,\ldots,0)$. Then define $m^*(X):=d_2^X$, which is called the \textbf{Euler multiplicity} on $\A^H$.
\label{Eulerest}
\end{define}

Now let us introduce two useful exact sequences:

\begin{define}[Euler sequence, Proposition 2.4, \cite{AD}]
Let $H \in \A$. 

(1)\,\,
For each $0 \le p \le \ell$, there is an exact sequence 
$$
0 \rightarrow 
\Omega^p(\A,m)[-1] \stackrel{\cdot \alpha_H}{\rightarrow }
\Omega^p(\A,m-\delta_H)
\stackrel{i_H^p}{\rightarrow }\Omega^p(\A^H,m^*)
$$
called the \textbf{Euler exact sequence}. Here the 
\textbf{Euler restriction map} $i=i^p=i_H^p$ is defined by, for $\alpha_H=x_1$ and for 
\begin{eqnarray*}
\omega &=&\sum_{1<i_2<\cdots <i_p \le \ell} 
(f_{i_2,\ldots,i_p}/Q(\A,m-\delta_H)) dx_1 dx_{i_2}\cdots dx_{i_p}\\
&\ &+
\sum_{1<i_1<\cdots <i_p \le \ell} 
(g_{i_1,i_2,\ldots,i_p}/Q(\A,m-\delta_H) )dx_{i_1} dx_{i_2}\cdots dx_{i_p}
\in \Omega^p(\A,m-\delta_H),
\end{eqnarray*}
as follows:
$$
i(\omega):=
\sum_{1<i_1<\cdots <i_p \le \ell} 
(\overline{g_{i_1,i_2,\ldots,i_p}/Q(\A,m-\delta_H))} dx_{i_1} dx_{i_2}\cdots dx_{i_p}, 
$$
where for $f \in S$, $\overline{f}$ denotes the image of $f$ in $\overline{S}=S/\alpha_H S$.

(2)\,\,
For each $0 \le p \le \ell$, there is an exact sequence 
$$
0 \rightarrow 
D^p(\A,m-\delta_H)[-1] \stackrel{\cdot \alpha_H}{\rightarrow }
D^p(\A,m)
\stackrel{\rho_H^p}{\rightarrow }D^p(\A^H,m^*)
$$
called the \textbf{Euler exact sequence}. Here the 
\textbf{Euler restriction map} $\rho=\rho^p=\rho_H^p$ is defined by 
$$
\rho(\theta)(\overline{f}_1,\ldots,\overline{f}_p):=
\overline{\theta(f_1,\ldots,f_p)}.
$$
\label{Euler}
\end{define}


Now let us introduce dual sequences of the Euler sequence, which was called the $B$-sequence introduced in \cite{A12} when $m\equiv 1$. Since we study several multiarrangements, we need its multiversion, but for that we need to modify $B$-sequences a bit, which we call the 
\textbf{$C$-sequences}.

\begin{theorem}[$C$-sequence]
Let $H \in \A$ and let $\delta_H:\A \rightarrow \{0,1\}$ be defined by $\delta_H(L):=\delta_{H,L}$.

(1)\,\,For each $0 \le p \le \ell$, there is an exact sequence 
$$
0 \rightarrow 
\Omega^p(\A,m-\delta_H) \stackrel{j}{\rightarrow }
\Omega^p(\A,m)
\stackrel{\res_H^p}{\rightarrow }(\Omega^{p-1}(\A^H,m^*) /\overline{C})[m(H)]
$$
called the \textbf{$C$-sequence}. Here 
$$
\overline{C}:=\overline{({\frac{Q(\A,m-\delta_H)}{\alpha_H^{m(H)-1}Q(\A^H,m^*)}})}
$$
is the \textbf{polynomial $C$}, $j$ is a canonical inclusion, and 
the 
\textbf{residue map} $\res:=\res_H^p$ is defined by, for $\alpha_H=x_1$ and for 
$$
\Omega^p(\A,m) \ni \omega
=\omega_0 \wedge \frac{dx_1}{x_1^{m(H)}}
+\omega_1
$$
with $\omega_1 \in (\Omega^{p-1}_{V'})_{(0)}\ (V'=\langle x_2,\ldots,x_\ell\rangle)$ and $\omega_1 \in (\Omega^p_{V'})_{(0)}$, 
$$
\res(\omega):=\overline{\omega_0}.
$$

(2)\,\,For each $0 \le p \le \ell$, there is an exact sequence 
$$
0 \rightarrow 
D^p(\A,m) \stackrel{j}{\rightarrow }
D^p(\A,m-\delta_H)
\stackrel{\partial_H^p}{\rightarrow }(D^{p-1}(\A^H,m^*)\overline{C})[-m(H)+1]
$$
called the \textbf{$C$-sequence}. Here 
$j$ is a canonical inclusion, and $\partial:=\partial_H^p=\partial^p$ is defined by, 
for $\theta \in D^p(\A,m-\delta_H)$, 
$$
\partial(\theta)(\overline{f}_2,\ldots,\overline{f}_p):=
\overline{(\displaystyle \frac{\theta(\alpha_H,f_2,\ldots,f_p)}{\alpha_H^{m(H)-1}})}.
$$
\label{Bseq}
\end{theorem}

\noindent
\textbf{Proof}.
Apply the identification (Proposition \ref{idenfify}) to the Euler sequences (Definition \ref{Euler}). We need to divide/muitiply the power of $\alpha_H$ to take modulo $\alpha_H$. \owari
\medskip

\begin{rem}
When $m \equiv 1$, $C$-sequences were just $B$-sequences, where 
$B$ is Terao's polynomial $B$ in \cite{T1}, and sequences were exactly introduced in \cite{A12}. 
\end{rem}

Let us recall a useful decomposition:

\begin{lemma}[Remark 3.3, \cite{A0}]
(1)\,\,
For any $H \in \A$ and each $0 \le p \le \ell$, there is a decomposition
$$
\Omega^p(\A) \simeq \Omega^{p-1}(\A) \wedge \frac{d \alpha_H}{\alpha_H}
\oplus \Omega^{p}(\A) \wedge \frac{d \alpha_H}{\alpha_H}.
$$
Moreover, 
$$
\Omega^0(\A,m)=S
$$
and 
$$
\Omega^\ell(\A,m) \simeq S \frac{1}{Q(\A,m)}.
$$

(2)\,\,
For any $H \in \A$ and each $1 \le p \le \ell-1$, let 
$$
D_H^p(\A):=\{\theta \in D^p(\A) \mid 
\theta(\alpha_H,f_2,\ldots,f_p)=0\ (\forall f_2,\ldots,f_p \in S)\}.
$$
Then 
there is a decomposition
$$
D^p(\A) = D^{p-1}_H(\A)\wedge \theta_E 
\oplus D^p_H(\A),
$$
where $\theta_E:=\sum_{i=1}^\ell x_i \partial_{x_i}$ is the \textbf{Euler derivation}. 
Moreover, 
$$
D^0(\A,m)=S
$$
and 
$$
D^\ell(\A,m) \simeq  S Q(\A,m).
$$

\label{decomp}
\end{lemma}

There is a canonical way to construct a multiarrangement from a given arrangement $\A$ and $H \in \A$ due to Ziegler as follows:

\begin{define}[\cite{Z}]
For each $0 \le p \le \ell$ and 
$H \in \A$, the \textbf{Ziegler restriction} $(\A^H,m^H)$ of $\A$ onto $H$ is defined by 
$$
m^H(X):=|\{L \in \A \setminus \{H\} \mid L \cap H=X\}|
$$
for $X \in \A^H$. Then there is the \textbf{Ziegler exact sequence} 
$$
0 \rightarrow (\Omega^p(\A)\wedge \frac{dx_1}{x_1})[-1]
\stackrel{\cdot \alpha_H}{\rightarrow }\Omega^p(\A)\wedge \frac{dx_1}{x_1} 
\stackrel{\pi^p_H}{\rightarrow }
\Omega^p(\A^H,m^H).
$$
Here we put $\alpha_H=x_1$ and $\pi^p_H$ is the restriction of $\res_H^{p+1}$ 
onto $\Omega^p(\A)\wedge \frac{dx_1}{x_1} $ which is called the 
\textbf{Ziegler restriction map}. Also for the derivation module, there is the \textbf{Ziegler 
exact sequence}
$$
0 \rightarrow D^p_H(\A)
\stackrel{\cdot \alpha_H}{\rightarrow }D_H^p(\A)
\stackrel{\pi^p_H}{\rightarrow }
D^p(\A^H,m^H),
$$
where $\pi_H^p:=\rho_H^p|_{D_H^p(\A)}$ is the 
\textbf{Ziegler restriction map}.
\label{Zrest}
\end{define}

Related to Ziegler restrictions, we have the following results.

\begin{theorem}[\cite{Z}]
Assume that $\A$ is free with $\exp(\A)=(1,d_2,\ldots,d_\ell)$. Then $\pi^p_H$ is surjective for all $p$, and $(
\A^H,m^H)$ is free with exponents $(d_2,\ldots,d_\ell)$ for any $H \in \A$
.
\label{Zrestfree}
\end{theorem}

\begin{theorem}[\cite{Y1}]
Let $\ell \ge 4$. 
Assume that $(\A^H,m^H)$ is free for some $H \in \A$ and $\A$ is locally free along $H$. Then 
$\A$ is free.
\label{Ycriterion}
\end{theorem}

\begin{theorem}[Projective dimensional surjection theorem, Theorem 3.2, \cite{A12}]
Let $H \in \A$ and $\A$ and $\A':=\A \setminus \{H\}$  be an essential arrangement. \begin{itemize}
    \item [(1)]
Assume that $\pd_S \Omega^p(\A,m)<\ell-2$ and $(\A,m)$ is locally free along $H$, i.e., 
$(\A_X,m_X)$ is free for all $0 \neq X \in L(\A^H)$. Then the Euler restriction map $i_H^p$ 
in Definition \ref{Euler} 
is surjective. In particular, the statement holds true if $(\A,m)$ is free.
\item [(2)]
Assume that $\pd_S \Omega^p(\A,m-\delta_H)<\ell-2$ and $(\A,m-\delta_H)$ is locally free along $H$. Then the residue map $\res_H^p$ 
in Definition \ref{Bseq} is surjective. 
In particular, the statement holds true if $(\A,m-\delta_H)$ is free.
\end{itemize}
\label{FST}
\end{theorem}

\begin{prop}[Proposition 4.14, \cite{OT}] 
Let $V=V_1 \times V_2$ and $(\A_i,m_i)$ be a multiarrangement in $V_i$. Let 
$\A=\A_1 \times \A_2$. Then 
$$
\Omega^p(\A,m_1+m_2)=\bigoplus_{i+j=p} \Omega^i(\A_1,m_1) \otimes_\K \Omega^j(\A_2,m_2).
$$
\label{product}
\end{prop}

\begin{prop}[Corollary 1.13, \cite{Eis}]
Let $M$ be a finitely generated $S$-graded module and let $\widetilde{M}$ be an associated sheaf on $\P^{\ell-1}$. 
If $\depth_S M \ge 2$, then 
$$
\Gamma_*(\widetilde{M}):=
\bigoplus_{k \in \Z}\Gamma(\widetilde{M}(k))=M.
$$
\label{Eisenbud}
\end{prop}

\begin{lemma}
Let $M$ be a finitely generated $S$-graded module and let $\widetilde{M}$ be an associated sheaf on $\P^{\ell-1}$. 
If $\pd_S M \le \ell-3$, then 
$$
H^1(\widetilde{M}(k))=0
$$
for all $k \in \Z$.
\label{H1zero}
\end{lemma}

\noindent
\textbf{Proof}.
Let us prove by induction on $0 \le k:=\pd_S$ that $H^i(\widetilde{M}(n))=0$ for all $n$ and 
all $1 \le i \le \ell-2-k$.
When $k=0$, there is nothing to show. 
By the same reason, there is nothing to show if $\ell-1 \le 1$. So we may assume that $k \ge 1 $ and 
$\ell-1 \ge 2$.  Let 
$$
0 \rightarrow \widetilde{K} \rightarrow \widetilde{F_0} \rightarrow \widetilde{M} \rightarrow 0
$$
be the first part of a minimal free resolution. Then $\pd_S \widetilde{K} =k-1$, thus the induction hypothesis shows that $H^i(\widetilde{K}(n))=0$
for all $n$ and $1 \le i \le \ell-2-k+1=\ell-k-1$. 
Thus the long exact sequence of the cohomology shows that 
$H^i(\widetilde{K}(n))=0$
for all $n$ and $1 \le i \le \ell-2-k$. So at least the first cohomology vanishes when $k\le  \ell-3$. \owari
\medskip

\begin{lemma}
Let $0 \neq \alpha \in V^*$ and $\overline{S}:=S/\alpha S$. Then for a finitely generated 
$\overline{S}$-graded module $M$, it holds that 
$$
\pd_{\overline{S}} M+1=\pd_S M.
$$
\label{pd+1}
\end{lemma}

\noindent
\textbf{Proof}.
For any regular sequence $\overline{a}_1,\ldots,\overline{a}_n \in \overline{S}$ for $M$ as an 
$\overline{S}$-module, the sequence $a_1,\ldots,a_n$ in $S$ form a regular sequence for $M$ as an $S$-module, and vice versa. It is trivial because the action of $S$ onto $M$ 
is through $\overline{S}$. 
Hence Auslander-Buchsbaum formula completes the proof.\owari
\medskip

The following can be proved by using the same argument as in Proposition 3.4 in \cite{AD2}.

\begin{prop}
Let $\A=\A' \cup \{H\}$, $m:\A \rightarrow \Z_{>0}$ be a multiplicity, and assume that $H$ is generic, i.e., 
$\codim_V H \cap \bigcap_{i=1}^{\ell-2} H_i=\ell-1$ for all distinct 
$H_1,\ldots,H_{\ell-2} \in \A'$.

(1)\,\,Then 
$i^p:\Omega^p(\A,m-\delta_H) \rightarrow \Omega^p(\A^H,m^*)$
is locally surjective, i.e., it is surjective at every localization except for the origin. 

(2)\,\, In the $C$-sequence in Definition \ref{Bseq}, both 
$\res^p$ and $\partial^p$ are locally surjective.
\label{genericsurjective}
\end{prop}

\noindent
\textbf{Proof}. 
First note that in each case the polynomial $C$ is just $1$ by genericity of $H$. So $\Omega^{p-1}(\A^H,m^*)/\overline{B}=\Omega^{p-1}(\A^H,m^*)$. 
Also, by the definition of $m^*$ in Definition \ref{Eulerest}, it holds that 
$$
m^*(H\cap L)=m(L).
$$
Now let $X \in L_{\ell-2}(\A^H)$. Then for (1) and (2), it suffices to show the surjectivity after 
localizing at $X$. Since $H$ is generic, we may assume that 
$X:x_1=\cdots=x_{\ell-1}=0$ and $H:x_1=0$. Let $m(x_i=0)=:m_i$. Then $(\A_X,m_X)$ is free with basis 
$$
\{\frac{dx_i}{x_i^{m_i}}\}_{i=1}^{\ell-1}, dx_\ell.
$$
By Proposition \ref{freewedge}, it holds that 
$$
\Omega^p(\A_X,m_X)=\wedge^p \langle \{\frac{dx_i}{x_i^{m_i}}\}_{i=1}^{\ell-1}, dx_\ell\rangle_S.
$$
Also, $(\A^H,m^*)$ is defined by $\prod_{i=2}^{\ell-1}x_i^{m_i}=0$. So the definitions of $i^p,\ 
\res^p$ and $\partial^p$ immediately show the local surjectivity. \owari
\medskip

%
\medskip

\section{Addition-deletion theorem for tameness}

Now let us prove Theorem \ref{addtame} by showing the following more general statement.

\begin{theorem}
Assume that $(\A,m-\delta_H)$ is tame, and both $(\A,m)$ and $(\A,m-\delta_H)$ are locally free along $H$. Then both $(\A,m)$ and
$(\A^H,m^*)
$ are tame. 
\label{addtamemulti}
\end{theorem}

\noindent
\textbf{Proof}.
It suffices to show that $\pd_S \Omega^p(\A,m) \le p$ for 
$1 \le p \le \ell-3$ since they are reflexive, thus 
$$
\pd_S \Omega^p(\A,m) \le \ell-2
$$
for all $p$. 
Recall the following 
Euler and $C$-sequences:
\begin{eqnarray*}
0 &\rightarrow& 
\Omega^p(\A,m) \stackrel{\cdot \alpha_H}{\rightarrow} \Omega^p(\A,m-\delta_H)
\stackrel{i^p}{\rightarrow} \Omega^p(\A^H,m^*),\\
0 &\rightarrow& \Omega^p(\A,m-\delta_H) \rightarrow \Omega^p(\A,m)
\rightarrow \Omega^{p-1}(\A^H,m^*)/C \rightarrow 0.
\end{eqnarray*}
The right exactness of the $C$-sequence follows from Theorem \ref{FST} and the local freeness 
of $(\A,m-\delta_H)$ along $H$ since $p<\ell-2$. 
We prove that $i^p$ is surjective, $\pd_S \Omega^p(\A,m) \le p$ and 
$\pd_S \Omega^p(\A^H,m^*) \le p+1$ for all $p \ge 1$ by induction on $1 \le p$, which are exactly what we have to show by Lemma \ref{pd+1}.
Assume that $p=1$.
First by the exactness of the $B$-sequence with the Ext-long exact sequence, we know that $\pd_S \Omega^1(\A,m) \le 1$ since $\pd_S \Omega^0(\A^H,m^*)=\pd_S \overline{S} =1$. Next 
look at the Euler exact sequence. Since $\A$ is locally free along $H$ and $\pd_S \Omega^1(\A,m) \le 1$ by the previous argument, Theorem \ref{FST} shows that $i^1$ is surjective. By the Ext-exact sequence of the Euler sequence, we know that $\pd_S \Omega^1(\A^H,m^*) \le 2$.

Now assume that up to $k \le p-1 <\ell-4,\ 0<p-1$, it holds that $i^k$ is surjective, 
$\pd_S \Omega^k(\A,m) \le k$ and 
$\pd_S \Omega^k(\A^H,m^*) \le k+1$. 
First by the exactness of the $C$-sequence of $\Omega^p(\A,m)$,  
the Ext-long exact sequence shows that 
\begin{eqnarray*}
0&=&\mbox{Ext}_S^k(\Omega^{p-1}(\A^H,m^*)/\overline{C},S) \rightarrow
\mbox{Ext}_S^k(\Omega^{p}(\A,m),S) \\
&\rightarrow&
\mbox{Ext}_S^k(\Omega^{p}(\A,m-\delta_H),S) =0
\end{eqnarray*}
for $k \ge p+1$  by induction hypothesis and the tameness of $(\A,m-\delta_H)$.
Hence $\mbox{Ext}_S^k(\Omega^{p}(\A,m),S)=0$ for $k \ge p+1$, implying that $\pd_S \Omega^p(\A,m) \le p$. Next 
look at the Euler exact sequence for $p$. Since $(\A,m)$ is locally free along $H$ and $\pd_S \Omega^p(\A,m) \le p$ by the previous argument, Theorem \ref{FST} shows that $i^p$ is surjective. Then the Ext-exact sequence of the Euler sequence shows that 
\begin{eqnarray*}
0&=&\mbox{Ext}_S^{k-1}(\Omega^{p}(\A,m),S) \rightarrow 
\mbox{Ext}_S^k(\Omega^{p}(\A^H,m^*),S)\\ &\rightarrow&
\mbox{Ext}_S^k(\Omega^{p}(\A,m-\delta_H),S)=0
\end{eqnarray*}
for $k \ge p+2$ since $\pd_S \Omega^p(\A,m) \le p$ by the previous argument and $(\A,m-\delta_H)$ is tame. 
Hence $\mbox{Ext}_S^k(\Omega^{p}(\A^H,m^*),S)=0$ for $k \ge p+2$, implying that $\pd_S \Omega^p(\A^H,m^*) \le p+1$, which completes the proof.\owari
\medskip

Let us show that a multi-generic arrangement is tame by using Theorem \ref{adddeltamemulti}.

\begin{theorem}
Let $\A$ be generic and irreducible. Then for any multiplicity $m$, the multiarrangement $(\A,m)$ is tame.
\label{multigenerictame}
\end{theorem}

\noindent
\textbf{Proof}.
Apply Therorem \ref{addtamemulti} to generic multiarrangements to complete the proof. For that purpose, we need to show that a generic multiarrangement is locally free, which is clear since every localization of such a 
multiarrangement at $0 \neq X$ is, after a change of coordinates, of the form 
$$
\prod_{i=1}^k x_i^{m_i}=0,
$$
which is clearly free. 
\owari
\medskip

Next let us show Theorem \ref{adddeltame}, which follows from the following multi-version.

\begin{theorem}
Assume that 
$(\A^H,m^*)$ is tame, and $(\A,m-\delta_H)$ is locally free along $H$. Then $(\A,m)$ is tame if 
$(\A,m-\delta_H)$ is tame.
\label{adddeltamemulti}
\end{theorem}

\noindent
\textbf{Proof}.
Since $(\A,m-\delta_H)$ is locally free and tame, Theorem \ref{FST} gives us the following exact sequence for $p<\ell-2$:
$$
0 \rightarrow \Omega^p(\A,m-\delta_H) \rightarrow \Omega^p(\A,m) 
\rightarrow \Omega^{p-1}(\A^H,m^*)/\overline{C} \rightarrow 0.
$$
Note that $\Omega^p(\A,m) \le p$ if $p=0,\ell-2,\ell-1,\ell$ by definitions and reflexivity of the logarithmic modules. Then applying the Ext-long exact sequence, 
the Ext-long exact sequence shows that 
\begin{eqnarray*}
0&=&\mbox{Ext}_S^k(\Omega^{p-1}(\A^H,m^*)/\overline{C},S) \rightarrow
\mbox{Ext}_S^k(\Omega^{p}(\A,m),S) \\
&\rightarrow&
\mbox{Ext}_S^k(\Omega^{p}(\A,m-\delta_H),S) =0
\end{eqnarray*}
for $k \ge p+1$ since both $(\A,m-\delta_H)$ and $(\A^H,m^*)$ are tame. 
Hence $\mbox{Ext}_S^k(\Omega^{p}(\A,m),S)=0$ for $k \ge p+1$, implying that $\pd_S \Omega^p(\A,m) \le p$, which 
completes the proof. \owari
\medskip




%

We can show Theorem \ref{resttame} in more general setup as follows:

\begin{theorem}
Assume that $(\A,m)$ and $(\A,m-\delta_H)$ are both tame, and 
$(\A,m)$ is locally free along $H$. Then $(\A^H,m^*)$ is tame.
\label{multitamerest}
\end{theorem}

\noindent
\textbf{Proof}.
Since $(\A,m)$ and $(\A,m-\delta_H)$ are tame and $(\A,m)$ is locally free along $H$, we have the right exact Euler sequence by Theorem \ref{FST}. Thus the Ext-computation shows that $\pd_S \Omega^p(\A^H,m^*) \le p+1$ for $p \le \ell-3$. Combining the reflexivity of $\Omega^p(\A^H,m^*)$, we can show that $(\A^H,m^*)$ is tame. \owari
\medskip

On the deletion theorem for tameness, the condition becomes a bit strong.

\begin{theorem}
Let $(\A,m)$ be tame, $H \in \A$ and assume that $(\A^H,m^*)$ is free, and $(\A,m)$ is locally free along $H$. 
Then $(\A,m-\delta_H)$ is tame.
\label{tamedeletion}
\end{theorem}

\noindent
\textbf{Proof}.
By the local freeness and Theorem \ref{FST}, we have
$$
0 \rightarrow \Omega^p(\A,m) \stackrel{\cdot \alpha_H}{\rightarrow} \Omega^p(\A,m-\delta_H) 
\rightarrow \Omega^p(\A^H,m^*) \rightarrow 0
$$
for $p=1,\ldots, \ell-3$. 
Now applying the Ext-exact and using $\pd_S \Omega^p(\A^H,m^*)=1$ by the freeness, we know that 
$(\A,m-\delta_H)$ is tame too by the same argument as in Theorem \ref{adddeltame}.\owari
\medskip

The following is immediate from Theorem \ref{adddeltame}.

\begin{cor}


Assume that $(\A,m)$ and $(\A^H,m^*)$ are free. Then $(\A,m-\delta_H)$ is tame.
\label{freetame}
\end{cor}

\begin{example}
Let $\A$ be the Edelman-Reiner arrangement in $\R^5$ 
introduced in \cite{ER} defined by 
$$
\prod_{i=1}^5 x_i \prod (x_1 \pm x_2 \pm x_3\pm x_4 \pm x_5)=0.
$$
It is free with exponents $(1,5,5,5,5)$. Let $H:x_1=0$. Then $\A^H$ 
is free with exponetns $(1,3,3,5)$. So Terao's deletion theorem shows that $\A':=
\A \setminus \{H\}$ is not free, but Corollary \ref{freetame} shows that $\A'$ is tame.
\end{example}

\section{Multiarrangements and tameness}


Since we know Ziegler's restriction theorem (Theorem \ref{Zrestfree}) and and Yoshinaga's criterion 
(Theorem \ref{Ycriterion}), it is natural to consider how tameness behaves when we take Ziegler restriction of tame 
arrangements, or when Ziegler restriction is tame. The former is nothing but Theorem \ref{Zresttame}.
\medskip

\noindent
\textbf{Proof of Theorem \ref{Zresttame}}.
By the local freeness, Definition \ref{Zrest} and Theorem \ref{Zrestfree}, we have the sheaf exact sequences
$$
0 \rightarrow \widetilde{\Omega^p(\A)\wedge 
\frac{d\alpha_H}{\alpha_H}}
\stackrel{\cdot \alpha_H}{\rightarrow }
\widetilde{\Omega^p(\A)\wedge 
\frac{d\alpha_H}{\alpha_H}}
\stackrel{\pi^H}{\rightarrow }
\widetilde{\Omega^p(\A^H,m^H)} \rightarrow 0
$$
for $0 \le p \le \ell-3$.
By Lemma \ref{pd+1}, it suffices to show that $\pd_{S} \Omega^p(\A^H,m^H) \le p+1$ for $1 \le p \le 
\ell-3$. 
First by Lemma \ref{decomp} and the tameness, we know that 
$$
\pd_S \Omega^p(\A)\wedge 
\frac{d\alpha_H}{\alpha_H} \le p
$$
for all $p$.
Thus $H^1(\widetilde{\Omega^p(\A)\wedge 
\frac{d\alpha_H}{\alpha_H}}(k))=0$ for all $k$ by Lemma \ref{H1zero}.
So by Lemma \ref{Eisenbud} and Theorem \ref{FST}, we have the right exact sequence 
$$
0 \rightarrow \Omega^p(\A)\wedge 
\frac{d\alpha_H}{\alpha_H}
\stackrel{\cdot \alpha_H}{\rightarrow }
\Omega^p(\A)\wedge 
\frac{d\alpha_H}{\alpha_H}
\stackrel{\pi^H}{\rightarrow }
\Omega^p(\A^H,m^H) \rightarrow 0.
$$
So the Ext-long exact sequence tells us that 
\begin{eqnarray*}
0&=&\mbox{Ext}_S^k (\Omega^p(\A)\wedge 
\frac{d\alpha_H}{\alpha_H},S) \rightarrow 
\mbox{Ext}_S^{k+1} (\Omega^p(\A^H,m^H),S) \\
&\rightarrow &
\mbox{Ext}_S^{k+1} (\Omega^p(\A)\wedge 
\frac{d\alpha_H}{\alpha_H},S) =0
\end{eqnarray*}
for $k \ge p+1$. Thus $\mbox{Ext}_S^{k} (\Omega^p(\A^H,m^H),S)=0$ for 
$k \ge p+2$, implying that $\pd_{\overline{S}} \Omega^p(\A^H,m^H) \le p$ for all $p$, which completes the proof. \owari
\medskip

In fact the proof of Theorem \ref{Zresttame} implies the following.

\begin{cor}
Let $1 \le p \le \ell-1$ and $0 < k < \ell-2$. Assume that $\A$ is locally free along $H$.
Assume that $\pd_S \Omega^p(\A) \le k$. 
Then $\pd_{\overline{S}} \Omega^i(\A^H,m^H)\le k$ for $i=p-1,p$. 


\label{corproof}
\end{cor}

We give an application of Theorem \ref{Zresttame}. Recall that for a multiarrangement $(\A,m)$, we can define its characteristic polynomial 
$$
\chi(\A,m;t)=\sum_{i=0}^\ell (-1)^{\ell-i} b_i(\A,m) t^i
$$ as in \cite{ATW}, but the relation between $\chi(\A;t)$ and $\chi(\A^H,m^H;t)$ have not yet been well-studied. For example, the coefficients $b_0(\A),\ldots,b_\ell(\A)$ of $\chi(\A;t)$ are the Betti numbers of the complements of $\A$ when $\K=\CC$, but there are no such meanings up to now for $\chi(\A,m;t)$. In fact, we even do not know whether 
$b_i(\A,m) \ge 0$ or not in general. The only known case is the following in \cite{A0}:

\begin{theorem}[Theorem 1.1, \cite{A0}]
Assume that $\A$ and $(\A^H,m^H)$ are both tame. Let 
$$
\chi_0(\A;t):=\chi(\A;t)/(t-1)=\sum_{i=0}^{\ell-1} b_i^0 (-1)^{\ell-1-i} t^{\ell-1-i} 
$$
and let 
$$
\chi(\A^H,m^H;t)=\sum_{i=0}^{\ell-1} \sigma_i (-1)^{\ell-1-i} t^{\ell-1-i} .
$$
Then $$
b_i^0 \ge \sigma_i \ge 0$$
for all $i$.
\label{tamebettiA0}
\end{theorem}

It was very difficult to check whether $(\A^H,m^H)$ is tame or not. However, now 
we can use Theorem \ref{Zresttame} to apply Theorem \ref{tamebettiA0} just by assuming that $\A$ is tame.

\begin{cor}
Assume that $\A$ is tame and $\A$ is locally free along $H \in \A$. 
Then in the notation of Theorem \ref{tamebettiA0}, 
$$
b_i^0 \ge \sigma_i \ge 0$$
for all $i$.
\label{tamebetti}
\end{cor}

\noindent
\textbf{Proof}. 
Apply Theorem \ref{tamebettiA0} by confirming tameness of $(\A^H,m^H)$ by Theorem \ref{Zresttame}.\owari
\medskip

Then how about the converse of 
Theorem \ref{Zresttame}? Namely, if $(\A^H,m^H)$ is tame, then can we say something about the tameness of $\A$? That is nothing but Theorem \ref{Ytame}. For the proof let us introduce one new result based on Theorem \ref{MS} in \cite{MS}.


\begin{theorem}
Let $(\A,m)$ be locally free, i.e., $(\A_X,m_X)$ is free for all $0 \neq X \in L(\A)$. 
If $\pd_S \Omega^1(\A,m)=1$, then $\wedge^p \Omega^1(\A,m)=\Omega^p(\A,m)$ and $\pd_S \Omega^p(\A,m)=p$ for 
$0 \le p \le \ell-2$. In particular, $(\A,m)$ is tame .
\label{MSmulti}
\end{theorem}

To prove Theorem \ref{MSmulti}, let us recall several results from \cite{MS}.

\begin{lemma}[\cite{MS}, Lemma 5.2]
Let $M$ be an $S$-graded module. If $M$ is locally free, then there is an exact sequence 
$$
0 \rightarrow M \rightarrow F_1 \rightarrow \cdots \rightarrow F_{\ell-2},
$$
where $F_i$ is a $S$-graded free module.
\label{sy}
\end{lemma}

\begin{prop}[\cite{L}, \cite{MS}, page 716]
Let $M$ be a locally free $S$-graded module with 
$\pd_S M=1$. Then for $1 \le p \le \ell-2$, it holds that $\pd_S \wedge^p M = p$.
\label{Lebelt}
\end{prop}

\noindent
\textbf{Proof of Theorem \ref{MSmulti}}. 
The proof is the same as that of Theorem \ref{MS}. Let us give a proof for the completeness.
Letting $M:=\Omega^1(\A,m)$, the assumptions on $(\A,m)$ show that all the assumptions in Proposition \ref{Lebelt} are satisfied. So $\pd_S \wedge^p 
\Omega^1(\A,m) = p$. Now compare $\wedge^p \Omega^1(\A,m)$ and $\Omega^p(\A,m)$. We know that 
$\wedge^p \Omega^1(\A,m) \subset \Omega^p(\A,m)$. Also, since $(\A,m)$ is locally free, 
the sheaf $E:=\widetilde{\Omega^1(\A,m)}$ on $\P^{\ell-1}$ is locally free.
By Proposition \ref{freewedge}, 
as sheaves, there is an isomorphism
$$
\widetilde{\wedge^p \Omega^1(\A,m)} \simeq \wedge^p E \rightarrow E^p:=\widetilde{\Omega^p (\A,m)}.
$$
Since $\pd_S \Omega^p(\A,m) \le \ell-2$, we know that 
$\mbox{depth}_S \Omega^p(\A,m) \ge 2$. So Proposition \ref{Eisenbud} shows that 
$\Gamma_*(\widetilde{\Omega^p(\A,m)})=\Omega^p(\A,m)$. Also, by taking the global sections, we know that 
$$
\Gamma_*(\widetilde{\wedge^p \Omega^1(\A,m)}) \simeq \Gamma_*(\wedge^p E) \simeq \Gamma_*(E^p) =\Omega^p(\A,m).
$$
By Proposition \ref{Lebelt}, we also know that $\mbox{depth}_S \wedge^p \Omega^1(\A,m) \ge 2$. So again Proposition \ref{Eisenbud} shows that $\wedge^p \Omega^1(\A,m) =\Gamma_*(\widetilde{\wedge^p \Omega^1(\A,m)})$, which shows that $\wedge^p \Omega^1(\A,m) =\Omega^p(\A,m)$. \owari
\medskip
%

\begin{rem}
Note that, contrary to Corollary 6.4 and the last statement in Theorem 6.5 in \cite{MS}, we cannot say anything on the projective dimension of $\Omega^{\ell-1}(\A,m)$ and $D^{\ell-1}(\A,m) $ since 
there are no decomposition like Proposition \ref{decomp} 
for a multiarrangement. 
\end{rem}

Moreover let us recall the following.

\begin{theorem}[Theorem 2.3, \cite{Y1}]
Let $E$ be a reflexive sheaf which is locally free except for at most finitely many points 
on $\P^\ell$. If $\ell \ge 3$, then $H^1(E(k))=0$ for all $k\ll 0$. 
\label{Yzero2}
\end{theorem}

\noindent
\textbf{Proof of Theorem \ref{Ytame}}.
(1)\,\,
This is trivial by Theorem \ref{adddeltame} since all multiarrangements in $\K^3$ are tame.

(2)\,\,
First let $p=1$.
By local freeness along $H$, we know that $(\A^H,m^H)$ is locally free by 
Theorem \ref{Zrestfree}, and there is a sheaf exact sequence 
$$
0 \rightarrow \widetilde{\Omega^1(\A)\wedge 
\frac{d\alpha_H}{\alpha_H}}
\stackrel{\cdot \alpha_H}{\rightarrow }
\widetilde{\Omega^1(\A)\wedge 
\frac{d\alpha_H}{\alpha_H}}
\stackrel{\pi^1_H}{\rightarrow }
\widetilde{\Omega^1(\A^H,m^H)} \rightarrow 0.
$$
Since $\pd_{\overline{S}}\Omega^1(\A^H,m^H) \le 1 <\ell-3$, it holds that 
$H^1(\widetilde{\Omega^1(\A^H,m^H)}(k))=0$
for all $k$ by Proposition \ref{H1zero}. So the cohomology long exact sequence shows that 
$$
H^1(\widetilde{\Omega^1(\A)\wedge 
\frac{d\alpha_H}{\alpha_H}}(k-1))
\stackrel{\cdot \alpha_H}{\rightarrow }
H^1(\widetilde{\Omega^1(\A)\wedge 
\frac{d\alpha_H}{\alpha_H}}(k))
$$
is surjective. Since $\A$ is locally free along $H$, it is easy to see that $\Omega^1(\A)$ has finitely many points at which it is not free. Thus  Theorem \ref{Yzero2} shows that 
$H^1(\widetilde{\Omega^1(\A)\wedge 
\frac{d\alpha_H}{\alpha_H}}(k))=0$ for all 
$k \in \Z$, implying that 
$$
\Gamma_*(\widetilde{\Omega^1(\A)\wedge 
\frac{d\alpha_H}{\alpha_H}})
\stackrel{\pi}{\rightarrow}
\Gamma_*(\widetilde{\Omega^1(\A^H,m^H)})
$$
is surjective. Since $\Omega^1(\A)$ is reflexive, the decomposition in Lemma \ref{decomp} shows that $\Omega^1(\A)\wedge 
\frac{d\alpha_H}{\alpha_H}$ is reflexive too. Thus Proposition \ref{Eisenbud} shows that 
$\Gamma_*((\widetilde{\Omega^1(\A)\wedge 
\frac{d\alpha_H}{\alpha_H}})=\Omega^1(\A)\wedge 
\frac{d\alpha_H}{\alpha_H}$. So again with Proposition \ref{Eisenbud}, there is an exact sequence 
$$
0 \rightarrow
\Omega^1(\A)\wedge 
\frac{d\alpha_H}{\alpha_H}
\stackrel{\cdot \alpha_H}{\rightarrow}
\Omega^1(\A)\wedge 
\frac{d\alpha_H}{\alpha_H}
\stackrel{\pi^1_H}{\rightarrow }
\Omega^1(\A^H,m^H)\rightarrow 
0.
$$
Taking the Ext-exact sequence combined with the tameness of $(\A^H,m^H)$, we know  that 
$$
M^k:=\mbox{Ext}_S^{k}(\Omega^1(\A)\wedge 
\frac{d\alpha_H}{\alpha_H},S)
\stackrel{\cdot \alpha_H}{\rightarrow}
\mbox{Ext}_S^{k}(\Omega^1(\A)\wedge 
\frac{d\alpha_H}{\alpha_H},S)=M^k
$$
is surjective for $k\ge 2$. Since $M^k$ 
is finitely generated $S$-graded module, and $(\alpha_H)M^k=M=(x_1,\ldots,x_\ell)M^k$, Nakayama's lemma 
shows that $M^k=0$ for $k \ge 2$. So Lemma \ref{decomp} shows that 
$\pd_S \Omega^1(\A) \le 1$. 


By Theorem \ref{MSmulti} and local freeness of $(\A^H,m^H)$, we have $\wedge^p \Omega^1(\A^H,m^H)=\Omega^p(\A^H,m^H)$ for 
$1 \le p \le \ell-3$. Also, since $\pi_H^1$ is surjective, 
for a set of generators $\omega_1,\ldots,\omega_s $ for $\Omega^1(\A^H,m^H)$, 
there are $\psi_i:=\eta_i \wedge \frac{d\alpha_H}{\alpha_H}+\zeta_i \in \Omega^1(\A)$ for $i=1,\ldots,s$,
$\eta_i \in \frac{1}{Q(\A)}\Omega^0_{V'},\zeta_i \in \frac{1}{Q(\A)}\Omega^1_{V'}$ such that 
$\pi(\psi_i \wedge (d\alpha_H/\alpha_H))=\omega_i$ for all $i$, where $V'=\K[x_2,\ldots,x_\ell]$ after putting $\alpha_H=x_1$. By Theorem \ref{MSmulti}, the set $$\{ 
\omega_{i_1} \wedge \cdots \wedge \omega_{i_p}\}_{1 \le i_1 <\cdots <i_p \le \ell}
$$ generates $\Omega^p(\A^H,m^H)$.  
Then the logarithmic differential form $$
\psi_{i_1} \wedge \cdots \wedge \psi_{i_p} 
\wedge \displaystyle \frac{d\alpha_H}{\alpha_H} \in \Omega^p(\A) \wedge 
\displaystyle \frac{d\alpha_H}{\alpha_H} 
$$
will be sent to $\omega_{i_1}\wedge \cdots \wedge \omega_{i_p}$ by $\pi^p$. Hence 
$$
\Omega^p(\A) \wedge \frac{d\alpha_H}{\alpha_H} \stackrel{\pi^p_H}{\rightarrow }
\Omega^p(\A^H,m^H)$$
is surjective for $p=1,\ldots,\ell-3$.
Thus taking the Ext-exact sequence combined with the tameness of $(\A^H,m^H)$, we know that 
$$
M^{p,k}:=\mbox{Ext}_S^{k}(\Omega^p(\A)\wedge 
\frac{d\alpha_H}{\alpha_H},S)
\stackrel{\cdot \alpha_H}{\rightarrow}
\mbox{Ext}_S^{k}(\Omega^p(\A)\wedge 
\frac{d\alpha_H}{\alpha_H},S)=M^{p,k}
$$
is surjective for $k\ge p+1$, i.e., 
$\alpha_H M^{p,k}=(x_1,\ldots,x_\ell)M^{p,k}=M^{p,k}$. Since 
$M^{p,k}$ is finitely generated $S$-graded module, Nakayama' lemma shows that  $M^{p,k}=0$ for $k \ge p+1$. So $\pd_S \Omega^p(\A) \wedge \frac{d\alpha_H}{\alpha_H} \le p$ for $1 \le p \le \ell-3$. Thus Lemma \ref{decomp} shows that 
$\pd_S \Omega^p(\A) \le p$ for $1 \le p \le \ell-3$. Since $\Omega^p(\A)$ is reflexive, 
we know that $\pd_S \Omega^p(\A) \le \ell-2$ for $p \ge \ell-2$, which shows that $\A$ is tame.\owari
\medskip

An easy corollary of the proof of Theorem \ref{Ytame} is the following:

\begin{cor}
Assume that $\pd_S \Omega^1(\A^H,m^H)=1$. Then $\A$ is tame if and only if $\A$ is locally free along $H$.
\end{cor}

\noindent
\textbf{Proof}. The key of the proof of Theorem \ref{Ytame} is the fact that $\pd_S \Omega^1(\A^H,m^H)=1$ and local freeness of $(\A^H,m^H)$. So it is clear.\owari
\medskip

\section{Inductive tameness}

There have been several class of free arrangements in which the freeness is combinatorial, i.e., determined by $L(\A)$. Here let us introduce the similar class of tame arrangements. 

\begin{define}
Let $\A$ be an arrangement and $X$ be a subset of an ambient space. We say that $\A$ is \textbf{stairly 
locally free along $X$} if $\A_x$ is stairly free for all $0 \neq x \in X$. For the stair freeness, see \cite{A6}.
\label{SLF}
\end{define}

\begin{define}
Let $ \ell \in \Z_{\ge 1}$ and $\mathcal{IT}_\ell$ consist of arrangements defined as follows.
First if $\ell \le 3$, then all arrangements are in $\mathcal{IT}_\ell$.
Also for all $\ell$, the empty arrangement is in $\mathcal{IT}_\ell$. For $4 \le \ell$, an arrangement $\A$ in $\K^\ell$ belongs to $\mathcal{IT}_\ell$ if either 
\begin{itemize}
    \item [(1)]
there is $H \in \A$ such that 
$\A$ and $\A\setminus \{H\}$ are stairly locally free along $H$, and 
$\A \setminus \{H\} \in \mathcal{IT}_\ell$, or 
\item[(2)]
there is $L \in \A$ such that 
$\A\setminus \{L\}$ is stairly locally free along $L$, 
$\A \setminus \{H\} \in \mathcal{IT}_\ell$ and $\A^L \in \mathcal{IT}_{\ell-1}$.
\end{itemize}
Let 
$$
\mathcal{IT}:=\bigcup_{0 \le \ell} \mathcal{IT}_\ell
$$
and $\A \in \mathcal{IT}$ is called the \textbf{inductively tame arrangement}.
\label{IT}
\end{define}

\begin{theorem}
$\A \in \mathcal{IT}$ is tame, and the tameness of an inductively tame arrangement is combinatorial.
\label{ITcombin}
\end{theorem}

\noindent
\textbf{Proof}.
Note that the freeness of the stair free arrangement is combinatorial by \cite{A6}. 
Since empty arrangements are tame, Theorems \ref{addtame} 
and \ref{adddeltame} complete the proof. \owari
\medskip

Also the following arrangements are tame which is combinatorially determined.

\begin{theorem}
Let $\A=\A' \cup \{H\}$ and assume that $\A'$ is tame. If the tameness of $\A'$ is determined by $L(\A')$, and $H$ is generic, then the tameness of $\A$ is determined by $L(\A)$ too.
\label{comtameadd}
\end{theorem}

\noindent
\textbf{Proof}.
Since we can obtain local surjectivities of both Euler and $B$-sequences by Proposition \ref{genericsurjective}, we can apply Theorem \ref{addtame} to complete the proof. \owari
\medskip

\begin{example}
Let $\A=\A' \cup \{H_1,\ldots,H_n\}$ and assume that $\A'$ is starily free in the sense of \cite{A6}. If $H_i$ is generic with respect to 
$\A_i:=\A' \cup \{H_1,\ldots,H_i\}$ for $i=1,\ldots,n$, then $\A_i$ is combinatorially tame.

Explicitly, for example, let $\A$ be defined by 
$$
Q(\A)=xyzw(x-y)(x-z)(x-w)(y-z)(y-w)(z-w),
$$
which is free with exponents $(1,2,3,4)$. If we add 
$H_1:x+3y+5z+7z=0,H_2:x+3^2y+5^2z+10^2z=0$ and so on, then we obtain non-free 
combinatorially tame arrangements. 
\end{example}


\end{document}